\documentclass[12pt,letterpaper]{amsart}
\usepackage{amsmath,txfonts}
\usepackage{amssymb}
\usepackage{amsxtra}
\usepackage{amsthm, color}
\usepackage{txfonts}
\usepackage{graphicx}
\usepackage{times}
\usepackage{citeref}
\usepackage{tikz}
\usepackage{hyperref}
\usepackage{stmaryrd}
\usepackage[T3,T1]{fontenc}
\usepackage{pgfplots}

\usetikzlibrary{calc}

\usepackage{mathrsfs}
\usepackage{amsfonts}
\usepackage{amssymb}
\usepackage{ifthen}
\usepackage{graphicx}
\nonstopmode \numberwithin{equation}{section}

\newtheorem{thm}{Theorem}[section]
\newtheorem{lem}{Lemma}[section]
\newtheorem{cor}[thm]{Corollary}
\newtheorem{prop}[thm]{Proposition}

\newtheorem{step}{Step}[section]

\theoremstyle{definition}
\newtheorem{mlem}{Main lemma}[section]
\newtheorem{assertion}{Assertion}[section]
\newtheorem{cl}{Claim}[section]
\newtheorem{ca}{Case}[section]
\newtheorem{sca}{Subcase}[section]
\newtheorem{scl}{Subclaim}[section]
\newtheorem{conj}[thm]{Conjecture}
\newtheorem{fact}{Fact}[section]
\newtheorem{defn}[thm]{Definition}
\newtheorem{op}[thm]{Open Problem}
\newtheorem{prob}{Problem}[section]
\newtheorem{ques}{Question}[section]
\newtheorem{rem}[thm]{Remark}
\newtheorem{exam}[thm]{Example}

\numberwithin{equation}{section}

\newcounter {own}
\def\theown {\thesection       .\arabic{own}}

\newenvironment{pf}[1][]{%
 \vskip 3mm
 \noindent
 \ifthenelse{\equal{#1}{}}%
  {{\slshape Proof. }}%
  {{\slshape #1.} }%
 }%
{\qed\bigskip}

\newcounter{alphabet}
\renewcommand{\thealphabet}{\Alph{alphabet}}

\newenvironment{Thm}[1][]{\refstepcounter{alphabet}%
	\bigskip
	\noindent
	{\bf Theorem \thealphabet}%
	\ifthenelse{\equal{#1}{}}{}{ (#1)}%
	{\bf .} \itshape
}{\vskip 8pt}

\newcommand{\IB}{{\mathbb B}}

\def\be{\begin{equation}}
\def\ee{\end{equation}}

\newcommand{\ben}{\begin{enumerate}}
\newcommand{\een}{\end{enumerate}}

\newcommand{\blem}{\begin{lem}}
\newcommand{\elem}{\end{lem}}
\newcommand{\bthm}{\begin{thm}}
\newcommand{\ethm}{\end{thm}}
\newcommand{\bcor}{\begin{cor}}
\newcommand{\ecor}{\end{cor}}
\newcommand{\beg}{\begin{exam}}
\newcommand{\eeg}{\end{exam}}
\newcommand{\begs}{\begin{examples}}
\newcommand{\eegs}{\end{examples}}
\newcommand{\bdefe}{\begin{defn}}
\newcommand{\edefe}{\end{defn}}
\newcommand{\bprob}{\begin{prob}}
\newcommand{\eprob}{\end{prob}}
\newcommand{\bques}{\begin{ques}}
\newcommand{\eques}{\end{ques}}
\newcommand{\bei}{\begin{itemize}}
\newcommand{\eei}{\end{itemize}}
\newcommand{\bcon}{\begin{conj}}
\newcommand{\econ}{\end{conj}}
\newcommand{\bop}{\begin{op}}
\newcommand{\eop}{\end{op}}

\newcommand{\bas}{\begin{assertion}}
\newcommand{\eas}{\end{assertion}}

\newcommand{\bfa}{\begin{fact}}
\newcommand{\efa}{\end{fact}}

\newcommand{\bca}{\begin{ca}}
\newcommand{\eca}{\end{ca}}

\newcommand{\bst}{\begin{step}}
\newcommand{\est}{\end{step}}

\newcommand{\bsca}{\begin{sca}}
\newcommand{\esca}{\end{sca}}

\newcommand{\bcl}{\begin{cl}}
\newcommand{\ecl}{\end{cl}}

\newcommand{\bmlem}{\begin{mlem}}
\newcommand{\emlem}{\end{mlem}}

\newcommand{\bscl}{\begin{scl}}
\newcommand{\escl}{\end{scl}}

\newcommand{\bcons}{\begin{conjs}}
\newcommand{\econs}{\end{conjs}}

\newcommand{\bprop}{\begin{prop}}
\newcommand{\eprop}{\end{prop}}

\newcommand{\br}{\begin{rem}}
\newcommand{\er}{\end{rem}}
\newcommand{\brs}{\begin{rems}}
\newcommand{\ers}{\end{rems}}
\newcommand{\bo}{\begin{obser}}
\newcommand{\eo}{\end{obser}}
\newcommand{\bos}{\begin{obsers}}
\newcommand{\eos}{\end{obsers}}
\newcommand{\bpf}{\begin{pf}}
\newcommand{\epf}{\end{pf}}
\newcommand{\ba}{\begin{array}}
\newcommand{\ea}{\end{array}}
\newcommand{\beq}{\begin{eqnarray}}
\newcommand{\beqq}{\begin{eqnarray*}}
\newcommand{\eeq}{\end{eqnarray}}
\newcommand{\eeqq}{\end{eqnarray*}}

\newcounter{minutes}
\divide\time by 60
\newcounter{hours}
\multiply\time by 60 \addtocounter{minutes}{-\time}

\begin{document}

\bibliographystyle{amsplain}
\title []
{Some sharp Schwarz type estimates  and their applications in Banach spaces}

\def\thefootnote{}
\footnotetext{ \texttt{\tiny File:~\jobname .tex,
		printed: \number\day-\number\month-\number\year,
		\thehours.\ifnum\theminutes<10{0}\fi\theminutes}
} \makeatletter\def\thefootnote{\@arabic\c@footnote}\makeatother
\author{Shaolin Chen}
\address{S. L. Chen,    Center for Applied Mathematics of Guangxi, Guangxi Normal University,
Guilin, Guangxi 541004, People's Republic of China} \email{mathechen@126.com}

\author{Hidetaka Hamada}
\address{H. Hamada, Faculty of Science and Engineering, Kyushu Sangyo University,
3-1 Matsukadai 2-Chome, Higashi-ku, Fukuoka 813-8503, Japan.}
\email{hi.hamada01@gmail.com; h.hamada@ip.kyusan-u.ac.jp}

\author{Megha Kundathil}
\address{M. Kundathil, Department of Mathematics, National Institute of Technology Calicut,
	Kozhikode--673 601, India.}
\email{meghakundathil10@gmail.com}

\author[R. Vijayakumar]{Ramakrishnan Vijayakumar}
\address{
	R. Vijayakumar, Department of Mathematics, National Institute of Technology Calicut,
	Kozhikode--673 601, India.}
\email{vijayakumar@nitc.ac.in}



	
	\makeatletter\def\thefootnote{\@arabic\c@footnote}\makeatother

\begin{abstract}
The primary objective of this paper is to develop methodologies for investigating Schwarz type lemmas and to present their applications in Banach spaces.
First, we improve upon the main results obtained by
Rogosinski [Jber. Dtsch. Math.-Verein 44: 258--261, 1934]
and Chen et al. [J. Anal. Math. 152: 181-216, 2024].
Based on these sharp estimates, we then derive several
sharp boundary Schwarz type lemmas (also known as Hopf type lemmas) for holomorphic mappings in Banach spaces, as well as for solutions
to certain classes of elliptic partial differential equations on the Euclidean unit ball in $\mathbb{C}^n$ or on the unit disk in $\mathbb{C}$. Furthermore, we prove some sharp Schwarz type lemmas
 for holomorphic mappings that send a prescribed point to another prescribed point. Finally, these lemmas are
 applied to establish a sharp Minda type Schwarz inequality in Banach spaces and to provide a sharp refined bound on subballs of the unit ball.

\end{abstract}

\maketitle \pagestyle{myheadings}
\markboth{S. L. Chen, H. Hamada, M. Kundathil and R. Vijayakumar}{Some sharp Schwarz type estimates  and their applications in Banach spaces}

\section{Introduction and Preliminaries }\label{csw-sec1}
It is well known that the Schwarz lemma stands as one of the most fundamental results in complex analysis. Its significance mainly stems from four key aspects: first, it provides a sharp growth bound for fixed-point self-mappings of the unit disk centered at the origin (see \cite{Car,Sch,Wu}); second, its invariant formulation-the Schwarz-Pick lemma-incorporates hyperbolic geometry into the framework of complex analysis (see \cite{Ah,ERS,Y}); third, it entails a rigidity property, where equality holds if and only if the mapping is a rotation (see \cite{BKR,BK,L-T});
fourth, it serves as a cornerstone for the classification of automorphisms, the uniqueness of Riemann mappings, boundary regularity analysis, complex dynamics, and geometric function theory (see
\cite{BE,CG,CH22,VM2Chen,CLW,EJLS,ELRS,FK,GHK-JAM,H15AdvMath,HK2023ANSNP,Ni,Wu}).
While the statement of the Schwarz lemma appears deceptively straightforward, its implications are profound and far-reaching. This is essentially because the lemma encapsulates the non-positive curvature of the unit disk in the language of complex analysis. Over the past century, it has evolved into a pivotal theme permeating numerous branches of mathematical research.
This paper advances the study of the classical Schwarz lemmas  and their applications to holomorphic mappings in Banach spaces.
In order to state our main results, we need to recall some basic definitions and introduce some necessary terminologies.

Let $\mathbb{C}^{n}$ be the $n$-dimensional complex Euclidean space, which is naturally identified with the real Euclidean space $\mathbb{R}^{2n}$.
In particular, we write $\mathbb{C} := \mathbb{C}^{1}$ and $\mathbb{R} := \mathbb{R}^{1}$.
For any $z = (z_{1}, \ldots, z_{n})$ and $w = (w_{1}, \ldots, w_{n})$ in $\mathbb{C}^{n}$, the standard Hermitian scalar product and the corresponding Euclidean norm are defined, respectively, as
$$\langle z,w\rangle :=
\sum_{k=1}^nz_k\overline{w}_k~\mbox{and}~|z|:={\langle
z,z\rangle}^{1/2}.$$
 We denote by $\mathbb{B}^{n}$ the open unit ball in $\mathbb{C}^{n}$ centered at the origin. In particular, $\mathbb{D} := \mathbb{B}^{1}$ denotes the unit disk in the complex plane $\mathbb{C}$.

\subsection*{Holomorphic mappings in Banach spaces} For real or complex Banach spaces $X$ and $Y$ with norms  $\|\cdot\|_X$ and
 $\|\cdot\|_Y$, respectively, denote by $L(X,Y)$  the space of all continuous linear operators from $X$ into $Y$, equipped with the standard operator norm:

\[
\|A\|=\displaystyle \sup_{x \in X\setminus \{0\} } \frac{\|Ax\|_Y}{\|x\|_X},
\]
where $A \in L(X,Y)$. Thus, $L(X,Y)$ is a Banach space under this norm. Denote by $X^{*}$ the dual space of the real or complex Banach space $X$. For
$x\in X \setminus\{0\}$, let
$$
T(x)=\{l_x \in X^{*}:l_x(x)=\|x\|_X ~\text{and}~ \|l_x\|_{X^{*}}=1\}.
$$
Then the well-known Hahn-Banach theorem implies that $T(x)\neq \emptyset$.

Let $f$ be a mapping of a domain $\Omega \subset X$ into a real or complex Banach space $Y$, where $X$ is a complex Banach space. We say that $f$ is differentiable at $z \in \Omega$ if there exists a bounded real linear operator $Df(z): X \rightarrow Y$ such that
$$
\lim_{\|h\|_X\rightarrow 0^+}\frac{\|f(z+h)-f(z)-Df(z)h\|_Y}{\|h\|_X}=0.
$$
Here $Df(z)$ is called the Fr\'echet derivative of $f$ at $z$. If $Y$ is a complex Banach space and $Df(z)$ is bounded and  complex linear for each $z\in \Omega$, then $f$ is said to be holomorphic on $\Omega$ (cf. \cite{VM2Chen,Chu-2021,Chu,C-R,HK2023ANSNP,HKK-26}).

Throughout this paper,  we denote by $\mathscr{H}(\Omega, \mathbb{B}_Y)$ the set of all holomorphic mappings from $\Omega$ to $\mathbb{B}_Y$, where $\Omega$ is a domain in the complex Banach space  $X$ and $\mathbb{B}_Y$ is the open unit ball of the complex Banach space  $Y$. Furthermore, for a given point $z_0 \in \mathbb{D}$, we denote by $\phi_{z_{0}}$
  the M\"obius transformation of $\mathbb{D}$ onto itself that satisfies $\phi_{z_0}(0)=z_0$, defined
  explicitly as
\beqq\label{yjg-01}\phi_{z_{0}}(\zeta)=\frac{z_{0}-\zeta}{1-\overline{z_{0}}\zeta},~\zeta\in\mathbb{D}.\eeqq

\subsection*{The Schwarz type lemmas for holomorphic mappings in Banach Spaces fixing the origin}
The classical Schwarz lemma states that every holomorphic mapping $f$ from  $\mathbb{D}$ into itself with $f(0)=0$ satisfies
\be\label{eq-ji-1}
|f(z)|\leq|z|\ee
for
all $z\in\mathbb{D}$. Moreover, unless $f$ is a rotation, the strict
inequality $|f'(0)|<1$ holds, and $f$ maps each disk $\mathbb{D}_{r}:=\{z:~|z|<r<1\}$
into a strictly smaller disk. Later, Lindel\"of removed the assumption ``$f(0)=0$'' and
improved the classical Schwarz lemma of holomorphic mappings \eqref{eq-ji-1} to the following refined inequality:

\beqq\label{eq-Lin}|f(z)|\leq\frac{|f(0)|+|z|}{1+|f(0)||z|},~z\in\mathbb{D}.\eeqq


A sharper form of the classical Schwarz lemma was given by the well-known Rogosinski Lemma \cite{Ro}. Specifically, for any holomorphic mapping $f: \mathbb{D} \to \mathbb{D}$ satisfying $f(0)=0$, the following inequality holds:
\be\label{jm-1} |f(z)|\leq\frac{|f'(0)|+|z|}{1+|f'(0)||z|}|z|,\quad z\in\mathbb{D}.\ee

Inequality (\ref{jm-1}) was also established independently by Osserman \cite{VM2_Osserman}. Moreover, Osserman employed inequality (\ref{jm-1}) to derive a sharp boundary Schwarz lemma, which can be stated as follows:
Let $f$ be a holomorphic mapping of $\mathbb{D}$ into itself such that $f(0)=0$. If $f$ is holomorphic at $a\in\partial\mathbb{D}$ (or, more generally, if $f$ is differentiable at $a\in\partial\mathbb{D}$) and satisfies $|f(a)|=1$, then
\be\label{jm-y1}|f'(a)|\geq\frac{2}{1+|f'(0)|}.\ee
Indeed, inequality (\ref{jm-y1}) is an immediate consequence of Unkelbach \cite[Satz 2]{U}.


Very recently, Chen et al.  extended Lindel\"of's inequality 
to
the setting of Banach spaces in the following form.

\begin{Thm}{\rm (\cite[Theorem 2.2]{VM2Chen})}\label{thm-h1}
Suppose that $\mathbb{B}_{X}$ and $\mathbb{B}_{Y}$ are the unit balls of the complex
Banach spaces $X$ and $Y$, respectively. Let $f:\, \mathbb{B}_{X}\rightarrow
\overline{\mathbb{B}_{Y}}$ be a holomorphic mapping. Then
\beqq\|f(z)\|_{Y}\leq\frac{\|f(0)\|_{Y}+\|z\|_{X}}{1+\|f(0)\|_{Y}\|z\|_{X}} \ \ \mbox{for}~z\in
\mathbb{B}_{X}. \eeqq
This estimate is sharp
with equality possible for each value of $\| f(0)\|_Y$
and for each $z\in \mathbb{B}_X$.
\end{Thm}

In particular, if $f(0)=0$ in Theorem \ref{thm-h1}, then Theorem \ref{thm-h1} can be reduced to the following form. 
 

\begin{Thm}{\rm (\cite[Theorem 2.3]{VM2Chen})}\label{Thm-A}
{Let $f\in\mathscr{H}(\mathbb{B}_X, \mathbb{B}_Y)$ with $f(0)=0$,
where $\mathbb{B}_{X}$ and $\mathbb{B}_{Y}$ are the open unit balls of complex
Banach spaces $X$ and $Y$, respectively.}
 Then
\beqq
\|f(z)\|_{Y}\leq\frac{\|Df(0)\|+\|z\|_{X}}{1+\|Df(0)\|\|z\|_{X}}\|z\|_{X}
\eeqq
holds for all $z\in\mathbb{B}_{X}$.
\end{Thm}

By an application of Theorem \ref{Thm-A}, Chen et al. established the following generalization of  (\ref{jm-y1}).

\begin{Thm}{\rm (\cite[Theorem 2.5]{VM2Chen})}\label{Thm-3h}
{Let $f\in\mathscr{H}(\mathbb{B}_X, \mathbb{B}_Y)$ with $f(0)=0$,
where $\mathbb{B}_{X}$ and $\mathbb{B}_{Y}$ are the open unit balls of complex
Banach spaces $X$ and $Y$, respectively.}
  If $f$ is holomorphic at $a\in\partial \mathbb{B}_X$
$($or more generally, if $f$ is differentiable at $a\in\partial \mathbb{B}_X$$)$
with
$\|f(a)\|_{Y}=1$, then \beqq \|Df(a)a\|_Y\geq\frac{2}{1+\|Df(0)\|}.
\eeqq
This inequality is sharp
with equality possible for each value of $\| Df(0)\|$.
\end{Thm}

The first purpose of this paper is to improve Theorems \ref{Thm-A} and \ref{Thm-3h} as follows.

\begin{thm}\label{thm-0}
	{
	Let $f\in\mathscr{H}(\mathbb{B}_X, \mathbb{B}_Y)$ with $f(0)=0$,
where $\mathbb{B}_{X}$ and $\mathbb{B}_{Y}$ are the open unit balls of complex
Banach spaces $X$ and $Y$, respectively.} Then, for every $z\in \mathbb{B}_X\setminus \{ 0\}$, we have
	\beqq
		\|f(z)\|_{Y}\leq\frac{\|Df(0)\eta\|_Y+\gamma_{f,\eta}(\|z\|_{X})}{1+\gamma_{f,\eta}(\|z\|_{X})\|Df(0)\eta\|_Y}\|z\|_{X},
	\eeqq
	where $\eta=z/\|z\|_{X},$ $$\gamma_{f,\eta}(t)=\frac{t+\lambda_{f,\eta}(0)}{1+\lambda_{f,\eta}(0)t}t,
	\quad
	t\in [0,1)$$ and
\begin{align}
\label{eq-lambda}
\lambda_{f,\eta}(0)&=
\left\{\begin{array}{ll}
\frac{\| D^2f(0)(\eta^2)\|_Y}{2(1-\|Df(0)\eta\|_Y^2)}, & \text { if \ }\|Df(0)\eta\|_Y<1,\\
1, & \text { if \ } \|Df(0)\eta\|_Y=1.
\end{array}\right.
 \end{align}
	Moreover, this inequality is sharp with equality possible for each values of
	$\| Df(0)\eta\|_{Y}\in [0,1)$ and $\lambda_{f,\eta}(0)\in [0,1)$.
\end{thm}

We remark that Theorem \ref{thm-0} is also a generalization of the inequality appearing in the seventh line from the bottom on page 511 of \cite{VM2Mercer}.
Applying Theorem \ref{thm-0}, we  obtain a sharp boundary Schwarz lemma as follows.

\begin{thm}\label{thm-3f}
{
Let $f\in\mathscr{H}(\mathbb{B}_X, \mathbb{B}_Y)$ with $f(0)=0$,
where $\mathbb{B}_{X}$ and $\mathbb{B}_{Y}$ are the unit balls of complex
Banach spaces $X$ and $Y$, respectively.}
  If
 $f$ is holomorphic at $a\in\partial \mathbb{B}_X$
$($or more generally, if $f$ is differentiable at $a\in\partial \mathbb{B}_X$$)$
with
$\|f(a)\|_{Y}=1$, then
\beqq \|Df(a)a\|_Y\geq1+\frac{2}{1+\lambda_{f,a}(0)}\left(\frac{1-\|Df(0)a\|_{Y}}{1+\|Df(0)a\|_{Y}}\right),
\eeqq where $\lambda_{f,a}(0)$ is defined in \eqref{eq-lambda}.
This inequality is sharp
with equality possible for each values of $\| Df(0)a\|_{Y}\in [0,1)$ and $\lambda_{f,a}(0)\in [0,1)$.
\end{thm}

We remark that if $\lambda_{f,a}(0)<1$, then
$$1+\frac{2}{1+\lambda_{f,a}(0)}\left(\frac{1-\|Df(0)a\|_{Y}}{1+\|Df(0)a\|_{Y}}\right)
>1+\frac{2}{1+1}\left(\frac{1-\|Df(0)\|}{1+\|Df(0)\|}\right)
=\frac{2}{1+\|Df(0)\|}.$$
Hence Theorem \ref{thm-3f} is an improvement of Theorem C
and \eqref{jm-y1} if $\lambda_{f,a}(0)<1$.

If we replace $\mathbb{B}_{X}$ and $\mathbb{B}_{Y}$ by $\mathbb{B}^{n}$ and $\mathbb{B}^{m}$ in Theorem \ref{thm-3f}, respectively,
 we obtain the following result, where $n$ and $m$ are positive integers.

\begin{cor}
Let $f\in\mathscr{H}(\mathbb{B}^{n}, \mathbb{B}^{m})$ with $f(0)=0$, where $n$ and $m$ are positive integers.
 If
 $f$ is holomorphic at $a\in\partial \mathbb{B}^{n}$
$($or more generally, if $f$ is differentiable at $a\in\partial \mathbb{B}^{n}$$)$
with
$|f(a)|=1$, then
\be\label{eq-a-01} |Df(a)a|\geq1+\frac{2}{1+\lambda_{f,a}(0)}\left(\frac{1-|Df(0)a|}{1+|Df(0)a|}\right),
\ee where $\lambda_{f,a}(0)$ is defined in \eqref{eq-lambda}.
This inequality is sharp
with equality possible for each values of $| Df(0)a|\in [0,1)$ and $\lambda_{f,a}(0)\in [0,1)$.

In particular, when $n=m=1$, {\rm (\ref{eq-a-01})} simplifies to the following sharper form:
\beqq
|f'(a)|\geq1+\frac{2}{1+\frac{|f''(0)|}{2(1-|f'(0)|^{2})}}\left(\frac{1-|f'(0)|}{1+|f'(0)|}\right).
\eeqq
\end{cor}


\subsection*{Applications of Schwarz type lemmas in partial differential equation theory}
Based on the idea of Theorem \ref{thm-3f}, we will study the boundary Schwarz type lemmas (also known as the Hopf type lemmas) for solutions
to certain classes of elliptic partial differential equations on $\mathbb{B}^n$ or on $\mathbb{D}$. To present our next results regarding the boundary Schwarz type lemma, we first recall some basic definitions.

Denote by $C^{m}(\Omega)$  the set of all  $m$-times continuously differentiable complex-valued
functions from $\Omega$ into $\mathbb{C}$, where $\Omega$ is a subset of $\mathbb{C}^n$
and $m$ is a nonnegative integer.
In particular, $C^0(\Omega)$ will be denoted by $C(\Omega)$.
For $f\in C^{2}({\Omega})$, let $$\Delta f(z)=\sum_{j=1}^n
\left(\frac{\partial^{2} f(z)}{\partial x_j^{2}}+\frac{\partial^{2} f(z)}{\partial y_j^{2}}\right)
$$ be the Laplacian of $f$,
where $z=(z_1,\dots, z_n)\in\Omega$ and $z_j=x_j+iy_j$ for $j\in \{ 1,\dots,n\}$. Since
$$\frac{\partial f}{\partial z_j}=\frac{1}{2}\left(\frac{\partial f}{\partial x_j}-i\frac{\partial f}{\partial y_j}\right)~\mbox{and}~ \frac{\partial f}{\partial\overline{z_j}}=\frac{1}{2}\left(\frac{\partial f}{\partial x_j}+i\frac{\partial f}{\partial y_j}\right),$$
we see that the Laplace operator can  be rewritten in the following form
 $$\Delta f(z)=4\sum_{j=1}^n\frac{\partial^{2} f(z)}{\partial z_j\partial \overline{z_j}}.$$

For a function $\psi\in C(\partial\mathbb{B}^n)$, we denote by $P[\psi]$ the Dirichlet solution for Laplace's operator
$\Delta$ of $\psi$ over $\mathbb{B}^n$, that is $\Delta P[\psi]=0$ in $\mathbb{B}^n$ and $P[\psi]=\psi$ on $\partial\mathbb{B}^n$.
It is well-known that, for  $\psi\in C(\partial\mathbb{B}^n)$,

$$P[\psi](z)=\int_{\partial\mathbb{B}^n}\mathbf{P}(z,\zeta)\psi(\zeta)d\sigma(\zeta),~z\in\mathbb{B}^n,$$
where $\mathbf{P}(z,\zeta)=\frac{1-|z|^{2}}{|\zeta-z|^{2n}}$
is the Poisson kernel
and $\sigma$ is the rotation-invariant positive Borel measure on $\partial\mathbb{B}^n$
for which $\sigma(\partial \mathbb{B}^n)=1$.

{
For   $g_{1}\in C(\overline{\mathbb{B}^n})$ and $f\in C^{2}(\mathbb{B}^n)$, the  inhomogeneous equation
$\Delta f=g_{1}$ is called Poisson's equation.
The Poisson equation is the fundamental law governing the static equilibrium of potential fields generated by distributed sources. Originating from 19th-century gravitational theory, it is the natural generalization of Laplace's equation and serves as the principal mathematical model for phenomena in electrostatics, gravitation, heat conduction, fluid mechanics, and beyond. Its elliptic nature and clear physical meaning make it a central pillar of both theoretical and computational mathematical physics (see \cite{C-H}).
In the following, we establish the boundary Schwarz type lemma (or the Hopf type lemma) for Dirichlet solutions to the Poisson equation.}

{\begin{thm}\label{thm-jc-1}
Suppose that $\varphi$ is holomorphic in $\mathbb{B}^n$ and continuous on $\overline{\mathbb{B}^n}$,
with $\varphi(0)=0$
and $\sup_{\zeta\in\partial \mathbb{B}^n}|\varphi(\zeta)|\leq1$. For a given $g_{1}\in C(\overline{\mathbb{B}^n})$,
let $f\in C^{2}(\mathbb{B}^n)$ satisfy the following boundary value problem:
\be\label{eq-v0-0}
	\left\{\begin{array}{ll}
\Delta f=g_{1}, & \text { in }\mathbb{B}^n,\\
		f=\varphi, & \text { in } \partial\mathbb{B}^n.
	\end{array}\right.
	\ee
If there exists some $\eta\in\partial\mathbb{B}^n$ such that $\lim_{r\rightarrow1^{-}}|f(r\eta)|=1$, then
\be\label{eq-wep-01}\liminf_{r\rightarrow1^{-}}\frac{|f(\eta)-f(r\eta)|}{1-r}\geq1+\frac{2}{1+\lambda_{\varphi,\eta}(0)}\left(\frac{1-|D\varphi(0)\eta|}{1+|D\varphi(0)\eta|}\right)-\frac{\|g_{1}\|_{\infty}}{2n}.\ee
Moreover,
the inequality of {\rm(\ref{eq-wep-01})} is sharp in the sense that
for each value of the right hand side of {\rm(\ref{eq-wep-01})} such that it is positive,
there exist $\varphi\in \mathscr{H}(\mathbb{B}^n, \mathbb{D})\cap C(\overline{\mathbb{B}^n})$ with $\varphi(0)=0$ and $g_{1}\in C(\mathbb{B}^n)$
such that the equality in \eqref{eq-wep-01} holds.
\end{thm}}

Let $\varphi_{1}$, $\varphi_{2}\in C(\partial\mathbb{D})$, $g\in C(\overline{\mathbb{D}})$,
and $f\in C^{4}(\mathbb{D})$. Of particular interest to us is the
following inhomogeneous biharmonic equation in $\mathbb{D}$:
$$\Delta(\Delta f)=g,$$
and its following associated Dirichlet boundary value problem:
$$
	\left\{\begin{array}{ll}
		\frac{\partial f}{\partial\overline{z}}=\varphi_{1}, & \text { in } \partial\mathbb{D}, \\
		f=\varphi_{2}, & \text { in } \partial\mathbb{D}.
	\end{array}\right.
	$$
The inhomogeneous biharmonic equation arises in areas of continuum mechanics,
including linear elasticity theory and the solution of Stokes flows (cf. \cite{Hay,Khu,We}).
Most important applications of the theory of functions of one complex variable were
obtained in the plane theory of elasticity and in the approximate theory of plates
subject to normal loading (cf. \cite{Ha,La}).


We now state the boundary Schwarz type lemma (or the Hopf type lemma) for solutions to the inhomogeneous biharmonic equation as follows.

\begin{thm}\label{thm-chmv-1}
 Suppose that $\varphi_{2}$ is holomorphic in $\mathbb{D}$ and continuous on $\overline{\mathbb{D}}$,
with $\varphi_{2}(0)=0$
and $\sup_{\zeta\in\partial \mathbb{D}}|\varphi_{2}(\zeta)|\leq1$.
For given $g\in C(\overline{\mathbb{D}})$ and $\varphi_{1}\in C(\partial\mathbb{D})$,
let $f\in C^{4}(\mathbb{D})$ satisfy the following boundary value problem:
\be\label{eq-vl-1}
	\left\{\begin{array}{ll}
\Delta(\Delta f)=g, & \text { in }\mathbb{D},\\
		\frac{\partial f}{\partial\overline{z}}=\varphi_{1}, & \text { in } \partial\mathbb{D}, \\
		f=\varphi_{2}, & \text { in } \partial\mathbb{D}.
	\end{array}\right.
	\ee 
If there exists some $\eta\in\partial\mathbb{D}$ such that $\lim_{r\rightarrow1^{-}}|f(r\eta)|=1$, then
\be\label{eq-wep-1}\liminf_{r\rightarrow1^{-}}\frac{|f(\eta)-f(r\eta)|}{1-r}\geq1+\frac{2}{1+\frac{|P''[\varphi_{2}](0)|}
{2(1-|P'[\varphi_{2}](0)|^{2})}}\left(\frac{1-|P'[\varphi_{2}](0)|}{1+|P'[\varphi_{2}](0)|}\right)-2|\varphi_{1}(\eta)|.\ee
Moreover,
the inequality of {\rm(\ref{eq-wep-1})} is sharp in the sense that
for each value of the right hand side of {\rm(\ref{eq-wep-1})} such that it is positive,
there exist $\varphi_{1}\in C(\partial\mathbb{D})$ and $\varphi_2\in \mathscr{H}(\mathbb{D}, \mathbb{D})\cap C(\overline{\mathbb{D}})$ with $\varphi_2(0)=0$
such that the equality in \eqref{eq-wep-1} holds for every $g\in C(\overline{\mathbb{D}})$.
\end{thm}


\subsection*{The boundary Schwarz type lemmas for holomorphic mappings without fixing the origin in Banach spaces}
We now recall a version of the boundary Schwarz type lemma that does not assume $f(0)=0$, where
$f$ is a holomorphic self-mapping of $\mathbb{D}$.  For related discussions, we refer the reader to \cite{VM2Chen}.

\begin{Thm}{\rm (\cite[Theorems 1.1 and 2.6]{Kalaj-arX})}\label{Thm-D}
{
Let $f\in\mathscr{H}(\mathbb{D}, \mathbb{B}^n)$ and $f(1)\in \partial \mathbb{B}^n$.
If $f'(1)$ exists,
then we have
\begin{equation}
\label{eq-Schwarz-Euclidean}
| f'(1)|\geq
\frac{2\left(1-| f(0)| \right)^2}
{1-| f(0)|^2+| f'(0)|}.
\end{equation}
\begin{enumerate}
\item[a)]
If $| f'(1)|=1$ and $f(0)=0$,
then $f$ is an affine disk.
\item[b)]
If $n=1$ and $f(0)=0$,
then the inequality \eqref{eq-Schwarz-Euclidean} reduces
\begin{equation}
\label{eq-Schwarz-disk}
| f'(1)|\geq
\frac{2}
{1+| f'(0)|},
\end{equation}
which is an equality if and only if $f(z)\equiv e^{i\theta}z$
or
\begin{align*}
f(z)\equiv e^{i\theta}z\frac{z+c}{1+{c}z}
\end{align*}
for some $\theta \in [0,2\pi)$ and $c\in [0,1)$.
\end{enumerate}
}
\end{Thm}

We remark that Theorem \ref{Thm-D} (b) was proved by Unkelbach \cite[Satz 2]{U}.


As our second goal, we seek to improve and generalize Theorem \ref{Thm-D}, presenting it in the form given below.

\begin{thm}\label{thm-Kalaj}
Suppose that
 $\Omega$ is a balanced domain in a complex Banach space $X$ and $\mathbb{B}_H$ is unit ball of a complex Hilbert space
$H$ with inner product $\langle \cdot, \cdot \rangle$.
Let $f\in\mathscr{H}(\Omega, \mathbb{B}_H)$.
If $f$ is differentiable at $z=a\in \partial \Omega$ and $f(a)=b\in \partial \mathbb{B}_H$,
then
\begin{equation}
\label{eq-Schwarz-Hilbert-new}
\langle Df(a)a, b\rangle\geq
\frac{2\left| 1-\langle f(0),b\rangle \right|^2}
{1-\left|\langle f(0),b \rangle \right|^2+|\langle Df(0)a, b\rangle|}.
\end{equation}

Moreover, if $f(0)=0$, then we have
\begin{equation}
\label{eq-Schwarz-Hilbert2-new}
\langle Df(a)a, b\rangle
\geq
\frac{2}
{1+|\langle Df(0)a, b\rangle|},
\end{equation}
which is an equality if and only if
$f(\zeta a)\equiv \zeta b$ or
\begin{align*}
f(\zeta a)\equiv \zeta\frac{\zeta+c}{1+{c}\zeta}b
\end{align*}
for some $c\in [0,1)$.

In particular, if $\Omega=\mathbb{B}_{X}$ is the unit ball of $X$,
then the inequality
$(\ref{eq-Schwarz-Hilbert-new})$ 
is sharp
with equality possible for each values of $\langle f(0),b\rangle$ and $|\langle Df(0)a, b\rangle|$
with $| \langle Df(0)a, b\rangle| \leq 1-|\langle f(0),b\rangle|^2$.
\end{thm}

Replacing the specific point ``$1$''
 in Theorem \ref{Thm-D} with an arbitrary point $z_{0}\in\partial\mathbb{D}$, we obtain the following slightly more general result in Banach spaces.

\begin{thm}\label{thm-09-CHMV}
Suppose that $\mathbb{B}_{X}$ and $\mathbb{B}_{Y}$ are the unit balls of the complex
Banach spaces $X$ and $Y$, respectively. Let $f:\, \mathbb{B}_{X}\rightarrow
\mathbb{B}_{Y}$ be a holomorphic mapping. If $f$ is differentiable at $z_{0}\in \partial \mathbb{B}_{X}$ and $f(z_{0})\in \partial \mathbb{B}_{Y}$,
then $$\|Df(z_{0})z_{0}\|_{Y}\geq\frac{1-\|f(0)\|_{Y}}{1+\|f(0)\|_{Y}}.$$ This estimate is sharp.
\end{thm}

\subsection*{Schwarz type lemmas and applications in Banach spaces concerning two prescribed points}
Our next objective is to derive a variant of the Schwarz lemma for holomorphic mappings $f: \mathbb{D} \to \mathbb{D}$
that fix the origin, i.e., $f(0)=0$, and send a prescribed point $z_{0}$ to $w_{0}$.
More specifically, suppose $f: \mathbb{D} \to \mathbb{D}$ is holomorphic with $f(0) = 0$, and suppose
there exist $a \in \mathbb{D} \setminus \{0\}$ and $\delta \in (0,1)$ such that $|f(a)| \leq |a| \delta$.
In this setting, Minda \cite{VM2Minda} obtained a Schwarz type inequality at the origin, yielding
 the bound:
\be\label{jm-3}
|f'(0)| \leq \frac{|a| + \delta}{1 + |a| \delta}.
\ee
This estimate has subsequently been applied to problems involving the hyperbolic metric and covering maps of Riemann surfaces.

In the following, we generalize the inequality given in (\ref{jm-3}) and establish its validity within the framework of Banach spaces.

\begin{thm}\label{thm-1}
	Let $f \in \mathscr{H}(\mathbb{D},\mathbb{B}_Y)$ with $f(0)=0$ and such that $f(z_0)=w_0~ (z_0\neq 0)$,
where $\mathbb{B}_{Y}$ is the open unit ball in a complex Banach space $Y$.
Then, for  $z \in \mathbb{D}$, we have

\begin{equation}\label{VM2_eq1}
		\|f(z)\|_{Y} \leq  \frac{(\|w_0\|_{Y}/|z_0|)+|\phi_{z_0}(z)|}{1+(\|w_0\|_{Y}/|z_0|)|\phi_{z_0}(z)|} |z|.
	\end{equation}
	Moreover, this inequality is sharp.
\end{thm}

We remark that in the one-dimensional case $\mathbb{B}_Y=\mathbb{D}$, Theorem \ref{thm-1} follows from a result of Mercer \cite[Proposition]{VM2Mercer}. 
Moreover, for the case $\mathbb{B}_Y=\mathbb{D}$ in Theorem \ref{thm-1}, we remove the assumption ``$f(0)=0$'' and obtain the following result.



\begin{thm}\label{VM2_tthm1}
    Let $f \in \mathscr{H}(\mathbb{D},\mathbb{D})$ with $f(z_0)=w_0~ (z_0\neq 0)$. Then, for  $z \in \mathbb{D}$, we have
	\begin{equation}\label{VM2_eqnq1}
		|f(z)| \leq \frac{K(z)|z|+|f(0)|}{1+K(z)|z||f(0)|},
	\end{equation}
	where
	$$
	K(z)=  \frac{(|\phi_{f(0)}(w_0)|/|z_0|)+|\phi_{z_0}(z)|}{1+(|\phi_{f(0)}(w_0)|/|z_0|)|\phi_{z_0}(z)|}.
	$$
	Moreover, this inequality is sharp.
\end{thm}

We remark that if $f(0)=0$ in Theorem \ref{VM2_tthm1}, then Theorem \ref{VM2_tthm1} reduces to Theorem \ref{thm-1} for $\mathbb{B}_Y=\mathbb{D}$.

Using Theorem \ref{thm-1}, we obtain a generalization of Minda's Schwarz type inequality (\ref{jm-3}), stated as follows.

\begin{thm} \label{VM2_thm6}
{
Let $f\in\mathscr{H}(\mathbb{B}_X, \mathbb{B}_Y)$ with $f(0)=0$,
where $\mathbb{B}_{X}$ and $\mathbb{B}_{Y}$ are the unit balls of complex
Banach spaces $X$ and $Y$, respectively.}
   Suppose that there exist $a \in \mathbb{B}_X \setminus \{0\}$ and $\delta \in (0,1)$ such that $\|f(a)\|_Y \leq \delta\|a\|_X $. Then the following sharp inequality holds:
$$ \left\| Df(0)\frac{a}{\|a\|_{X}}\right\|_{Y} \leq \frac{\|a\|_{X}+\delta}{1+\delta\|a\|_{X}}.
$$
In particular, if $\mathbb{B}_X=\mathbb{D}$, then the following sharp inequality holds:
$$ \left\| Df(0)\right\| \leq \frac{|a|+\delta}{1+\delta|a|}.
$$

\end{thm}


In what follows, we estimate the value of $f$ at $z$ by that at some point on the real interval  
from $0$
to $z/\| z\|_X$. This is a generalization of Theorem \ref{Thm-A}.

\begin{thm}\label{thm-s-1}
	{
	Let $f\in\mathscr{H}(\mathbb{B}_X, \mathbb{B}_Y)$ with $f(0)=0$,
where $\mathbb{B}_{X}$ and $\mathbb{B}_{Y}$ are the unit balls of complex
Banach spaces $X$ and $Y$, respectively.} Then, for every $z \in \mathbb{B}_X\setminus\{0\}$ and $r\in(0,1)$, we have
	\begin{equation}\label{VM2_eq5}
		\|f(z)\|_{Y} \leq \left( \frac{\frac{\left\|f(\frac{rz}{\|z\|_{X}})\right\|_{Y}}{r}+|\phi_r(\|z\|_{X})|}{1+\frac{\left\|f(\frac{rz}{\|z\|_{X}})\right\|_{Y}}{r}|\phi_r(\|z\|_{X})|} \right)\|z\|_{X}.
	\end{equation}
	Moreover, 
	 the above inequality is sharp.
\end{thm}

\begin{rem}
By letting $r\rightarrow0^{+}$ in (\ref{VM2_eq5}),
we can derive the result in \cite[Theorem 2.3]{VM2Chen}, which reads as follows:
\beqq
\|f(z)\|_{Y}\leq\frac{\|Df(0)\|+\|z\|_{X}}{1+\|Df(0)\|\|z\|_{X}}\|z\|_{X}.
\eeqq
Thus, Theorem \ref{thm-s-1} is a direct generalization of Theorem \ref{Thm-A}.
\end{rem}




By using Theorem \ref{thm-s-1}, we get the following estimate.

\begin{cor}\label{Vm2_thm1}
	{
	Let $f\in\mathscr{H}(\mathbb{B}_X, \mathbb{B}_Y)$ with $f(0)=0$,
where $\mathbb{B}_{X}$ and $\mathbb{B}_{Y}$ are the unit balls of complex
Banach spaces $X$ and $Y$, respectively.} Then for every $z \in \mathbb{B}_X\setminus \{0\}$ and $0<r<1$, we have
	\begin{equation*}
		\|Df(0)\| \leq \frac{\omega+r^2}{r(1+\omega)},
	\end{equation*}
	where $\omega = \sup\limits_{\|z\|_{X}=r}
	\left\| f(z) \right\|_{Y}$.
\end{cor}

The proofs of Theorems  \ref{thm-0}$\sim$\ref{thm-s-1} and Corollary \ref{Vm2_thm1} will be given in Section \ref{S-2}.

\section{The proofs of the main results}\label{S-2}





\begin{lem}\label{M-1} 
Let $f$ be a holomorphic mapping of $\mathbb{D}$ into itself with $f(0)=0$.
Then
	 \be\label{jm-2}
		|f(z)|\leq\frac{|f'(0)|+\nu(z)}{1+\nu(z)|f'(0)|}|z|\leq\frac{|f'(0)|+|z|}{1+|f'(0)||z|}|z|,~z \in \mathbb{D},
	\ee
	where $$\nu(z)=\frac{|z|+\frac{|f''(0)|}{2(1-|f'(0)|^2)}}{1+\frac{|f''(0)|}{2(1-|f'(0)|^2)}|z|}|z|.$$
\end{lem}

\bpf The inequality $$|f(z)|\leq\frac{|f'(0)|+\nu(z)}{1+\nu(z)|f'(0)|}|z|,~z\in\mathbb{D},$$
follows from \cite[Proposition]{VM2Mercer}.

Next, we will prove $$\nu(z)\leq|z|$$ for $z\in\mathbb{D}$.

Let us recall the following classical result which we call the coefficient type
Schwarz-Pick lemma of holomorphic mappings (see \cite{Ne}): if $\Phi$ is a holomorphic
mapping of $\mathbb{D}$ into itself with $\Phi(z)=\sum_{n=0}^{\infty}a_{n}z^{n}$, then
for each $n\geq1$,
\be\label{CS-1}|a_{n}|\leq1-|a_{0}|^{2}.\ee
Applying (\ref{CS-1}) to $f(z)/z$, we have
 $$\frac{|f''(0)|}{2(1-|f'(0)|^2)}\leq1,$$
which implies that $$\nu(z)\leq|z|$$ and
\beqq
|f(z)|\leq\frac{|f'(0)|+\nu(z)}{1+\nu(z)|f'(0)|}|z|\leq\frac{|f'(0)|+|z|}{1+|f'(0)||z|}|z|
\eeqq
which implies that (\ref{jm-2}) holds. The proof of this lemma is complete.
\epf

\subsection*{Proof of Theorem \ref{thm-0}}
Let $z$ be any fixed point in $\mathbb{B}_X\backslash\{0\}$
and let $\eta=z/\|z\|_{X}\in\partial\mathbb{B}_{X}$.
Without loss of generality,
we assume that $f(z)\neq0$ and $\|Df(0)\eta\|_{Y}<1$.
For any fixed $b\in\partial\mathbb{B}_{Y}$, let
$$F(\zeta)=l_{b}(f(\zeta\eta)),~\zeta\in\mathbb{D},$$
where $l_{b}\in T(b).$ Then $F$ is a holomorphic mapping from $\mathbb{D}$ into itself
with $F(0)=0$.
Applying Lemma \ref{M-1} to $F$, we get

\be\label{jm-20}
		|F(\zeta)|\leq \frac{|F'(0)|+\nu(\zeta)}{1+\nu(\zeta)|F'(0)|}|\zeta|,~\zeta \in \mathbb{D},
	\ee
	where $$\nu(\zeta)=\frac{|\zeta|+\frac{|F''(0)|}{2(1-|F'(0)|^2)}}{1+\frac{|F''(0)|}{2(1-|F'(0)|^2)}|\zeta|}|\zeta|.$$

Since $F'(\zeta)=l_{b}(Df(\zeta\eta)\eta)$ and $F''(\zeta)=l_{b}\left(D\big(Df(\zeta\eta)\eta\big)\eta\right)$,
we see that
\be\label{jm-21}
|F'(0)|=|l_{b}(Df(0)\eta)|\leq\|l_{b}\|_{Y^{\ast}}\|Df(0)\eta\|_{Y}=\|Df(0)\eta\|_{Y}<1
\ee
and
\be\label{jm-22}
|F''(0)|=\left|l_{b}\left(D^2f(0)(\eta^2)\right)\right|\leq\|l_{b}\|_{Y^{\ast}}\left\|D^2f(0)(\eta^2)\right\|_{Y}
=\left\|D^2f(0)(\eta^2)\right\|_{Y}.
\ee
Since the function
\be\label{hyu-1} 
Q(x)= \frac{x + c}{1 + cx} \ee is increasing with respect to $ x\in[0,\infty)$ for a fixed $ c \in [0,1)$, it follows from \eqref{jm-20}, (\ref{jm-21}) and (\ref{jm-22}) that

\beqq
\nu(\zeta)=\frac{|\zeta|+\frac{|F''(0)|}{2(1-|F'(0)|^2)}}{1+\frac{|F''(0)|}{2(1-|F'(0)|^2)}|\zeta|}|\zeta|
\leq\gamma_{f, \eta}(|\zeta|)
\eeqq
and
\beq\label{jm-23}
|l_{b}(f(\zeta\eta))|&=&|F(\zeta)|\leq \frac{\|Df(0)\eta\|_{Y}+\nu(\zeta)}{1+\nu(\zeta)\|Df(0)\eta\|_{Y}}|\zeta|\\ \nonumber
&\leq&\frac{\|Df(0)\eta\|_{Y}+\gamma_{f,\eta}(|\zeta|)}{1+\gamma_{f, \eta}(|\zeta|)\|Df(0)\eta\|_{Y}}|\zeta|.
\eeq
Finally, by letting $\zeta=\|z\|_{X}$ and $b=f(z)/\|f(z)\|_{Y}$ in (\ref{jm-23}), we get
\beqq
\|f(z)\|_{Y}\leq\frac{\|Df(0)\eta\|_{Y}+\gamma_{f,\eta}(\|z\|_{X})}{1+\gamma_{f,\eta}(\|z\|_{X})\|Df(0)\eta\|_{Y}}\|z\|_{X}.
\eeqq

Next, we prove the sharpness part.
Fix $\alpha,\beta\in[0,1)$ and $w\in\partial\IB_Y$.
Let $z_1\in\IB_X\setminus\{0\}$ and let $l_{z_1}\in T(z_1)$. 	
For $z\in\mathbb{B}_X$, let
\beqq
\label{eq-G}
f(z)=-F_{1}(-l_{z_{1}}(z))w,
\eeqq
where
\begin{align}
\label{def-F1}
F_{1}(\xi)&=\xi\phi_{\alpha}(\xi\phi_{\beta}(\xi)), \quad \xi\in\mathbb{D}.
\end{align}
Then $f$ is a holomorphic mapping of $\mathbb{B}_X$ into $\mathbb{B}_Y$
with \(f(0)=0\).
	
	Elementary computations show that for $\zeta\in \overline{\mathbb{D}}$
	and $\eta_{1}=z_{1}/\|z_{1}\|_{X}$,
\begin{align*}
f(\zeta \eta_1)&=-F_1(-\zeta)w
\\
&=\zeta\frac{\alpha+\zeta\phi_{\beta}(-\zeta)}{1+\alpha\zeta\phi_{\beta}(-\zeta)}w
\\
&=\zeta\frac{\alpha+\zeta\frac{\beta+\zeta}{1+\beta\zeta}}{1+\alpha\zeta\frac{\beta+\zeta}{1+\beta\zeta}}w,
\end{align*}
	\beq\label{eq-wed-1}
		\|f(z_1)\|_{Y}
		&=&\| f(\|z_1\|_X \eta_1)\|_{Y}\\ \nonumber
		&=&\|z_1\|_{X}
		\frac{\alpha+\|z_1\|_{X}\left(\frac{\beta+\|z_1\|_{X}}{1+\beta\|z_1\|_{X}}\right)}
		{1+\alpha\|z_1\|_{X}\left(\frac{\beta+\|z_1\|_{X}}{1+\beta\|z_1\|_{X}}\right)},
	\eeq
\be\label{yh-o1}Df(\zeta\eta_{1})\eta_{1}=F'_{1}(-\zeta)w,\ee
and
\be\label{yh-o2}D^2f(\zeta\eta_{1})(\eta_{1}^2)=-F''_{1}(-\zeta)w.\ee
It follows from (\ref{yh-o1}) and (\ref{yh-o2}) that
	$$Df(0)\eta_{1}
	=F'_{1}(0)w$$
and
$$D^2f(0)(\eta_{1}^2)=-F''_{1}(0)w,$$
which gives that
\beqq
|F'_{1}(0)|=\left\|Df(0)\eta_{1}\right\|_{Y}
\eeqq
and
\beqq
|F''_{1}(0)|=\|D^2f(0)(\eta_{1}^2)\|_{Y}.
\eeqq
Since $\alpha=|F'_{1}(0)|=\| Df(0)\eta_1\|_Y$ and
\beqq
\beta=\frac{|F''_{1}(0)|}{2(1-\alpha^{2})}=\lambda_{f,\eta_1}(0)
\eeqq
	hold,
	it follows from (\ref{eq-wed-1}) that
	
	\[
	\|f(z_1)\|_{Y}
	=\frac{\| Df(0)\eta_1\|_Y+\gamma_{f, \eta_1}(\|z_1\|_X)}
	{1+\gamma_{f,\eta_1}(\|z_1\|_X)\| Df(0)\eta_1\|_Y}
	\|z_1\|_X.
	\]
The proof of this theorem is complete.
\qed

\subsection*{Proof of Theorem \ref{thm-3f}}
By the triangle inequality, we have
\be\label{PV- Eq9} \left\|\frac{f(z)-f(a)}{\|z\|_{X} -
\|a\|_{X}}\right\|_{Y} \geq \frac{1-\|f(z)\|_{Y}}{1-\|z\|_{X}}\ee
for $z\in \mathbb{B}_X\setminus\{0\}$. Let $\eta=z/\|z\|_{X}$.
 It
follows from Theorem \ref{thm-0} that

\beqq\label{PV- Eq 10} \frac{1-\|f(z)\|_{Y}}{1-\|z\|_{X}} \geq
\frac{1-\frac{\left(\|Df(0)\eta\|_{Y}+\gamma_{f, \eta}(\|z\|_{X})\right)}{1+\gamma_{f,\eta}(\|z\|_{X})\|Df(0)\eta\|_{Y}}\|z\|_{X}}{1-\|z\|_{X}},
 \eeqq
 which, together with (\ref{PV- Eq9}), implies that

\be \label{PV- Eq 11}
 \liminf\limits_{r\rightarrow 1^{-}} \left\|\frac{f(ra)-f(a)}
 {\|ra\|_{X} - \|a\|_{X}}\right\|_{Y} \geq \liminf\limits_{r\rightarrow 1^{-}}\psi(r),
\ee
and
$$\psi(r)=\frac{1-\frac{\left(\|Df(0)a\|_{Y}+\gamma_{f,a}(\|ra\|_{X})\right)}{1+\gamma_{f,a}(\|ra\|_{X})\|Df(0)a\|_{Y}}\|ra\|_{X}}{1-\|ra\|_{X}}.$$
Since $$\lim_{r\rightarrow1^{-}}\gamma_{f,a}(r)=1,$$ by L'Hospital's Rule,
we see that
\beq\label{nx-01}
\liminf\limits_{r\rightarrow 1^{-}}\psi(r)
 &=&\lim_{r\rightarrow 1^{-}}\frac{1-\frac{\left(\|Df(0)a\|_{Y}+\gamma_{f,a}(r)\right)}{1+\gamma_{f,a}(r)\|Df(0)a\|_{Y}}r}{1-r}\\ \nonumber
 &=& \lim_{r\rightarrow 1^{-}}\frac{r\gamma_{f,a}'(r)\left(1-\|Df(0)a\|_{Y}^{2}\right)}{\left(1+\|Df(0)a\|_{Y}\gamma_{f,a}(r)\right)^{2}}
+
\lim_{r\rightarrow 1^{-}}\frac{\|Df(0)a\|_{Y}+\gamma_{f,a}(r)}{1+\|Df(0)a\|_{Y}\gamma_{f,a}(r)}\\ \nonumber
&=&1+\frac{2}{1+\lambda_{f,a}(0)}\frac{1-\|Df(0)a\|_{Y}}{1+\|Df(0)a\|_{Y}}.
\eeq
Since the radial derivative

\be\label{nx-02}
Df(a)a=\lim_{r\rightarrow 1^{-}}\frac{f(ra)-f(a)}{r-1}
\ee
exists, the desired result follows from 
(\ref{PV- Eq 11}), (\ref{nx-01}) and (\ref{nx-02}).

The sharpness part follows from the mapping $f(z)=-F_{1}(-l_{a}(z))w$,
where $w\in \partial \mathbb{B}_Y$ and $F_1$ is defined in \eqref{eq-G}.
The proof of this theorem is complete.
\qed


\subsection*{Proof of Theorem \ref{thm-jc-1}}
{It is well known that   all solutions to the Poisson
equation (\ref{eq-v0-0}) are given by (see \cite[p.118-120]{Ho})
\beqq
f(z)=P[\varphi](z)+\mathcal{G}[g_{1}](z),~z\in\mathbb{B}^n,
\eeqq where
$$\mathcal{G}[g_{1}](z)=\int_{\mathbb{B}^n}g_{1}(w)G(z,w)dV(w).$$
Here, $G(z,w)$ is the Green function and $dV(w)$ is the Lebesgue measure in $\mathbb{B}^n$.
It follows from
Theorem \ref{thm-0}  that
\beq\label{jut-1}
\left|f(r\eta)\right|&\leq&\frac{|D\varphi(0)\eta|+\gamma_{\varphi,\eta}(r)}{1+\gamma_{\varphi,\eta}(r)|D\varphi(0)\eta|}r
+\left|\mathcal{G}[g_{1}](r\eta)\right|,
\quad
r\in (0,1).
\eeq
By the uniqueness of solutions to the boundary value problem
\eqref{eq-v0-0} for $\varphi\equiv 0$ and $g_1\equiv 1$,
we have
\begin{align*}
\mathcal{G}[1](z)&\equiv \frac{1}{4n}(|z|^2-1),
\end{align*}
which combined with the fact that the Green function is non-positive
implies that
\beq\label{jut-2}
\left|\mathcal{G}[g_{1}](z)\right|&\leq&-{\|g_{1}\|_{\infty}}
\int_{\mathbb{B}^n}G(z,w)dV(w)\\ \nonumber
&=&\frac{\|g_{1}\|_{\infty}}{4n}(1-|z|^{2}),
\eeq
where $$\|g_{1}\|_{\infty}=\sup_{z\in\mathbb{B}^n}|g(z)|.$$
From (\ref{jut-1}) and (\ref{jut-2}), we obtain

\beqq
\liminf_{r\rightarrow1^{-}}\frac{|f(\eta)-f(r\eta)|}{1-r}&\geq&
\lim_{r\rightarrow1^{-}}\frac{1-\frac{|D\varphi(0)\eta|+\gamma_{\varphi,\eta}(r)}{1+\gamma_{\varphi,\eta}(r)|D\varphi(0)\eta|}r}{1-r}-
\frac{1}{4n}\|g_{1}\|_{\infty}\lim_{r\rightarrow1^{-}}\frac{1-r^{2}}{1-r}\\
&=&1+\frac{2}{1+\lambda_{\varphi,\eta}(0)}\frac{1-|D\varphi(0)\eta|}{1+|D\varphi(0)\eta|}-
\frac{\|g_{1}\|_{\infty}}{2n}.
\eeqq

Next, we come to prove the sharpness part.  For $z\in\overline{\mathbb{B}^n}$, let
\beqq
\label{eq-G2}
\varphi(z)=-F_{1}(-l_{\eta}(z)),
\eeqq
where $l_{\eta}\in T(\eta)$ and $F_{1}$ is defined in \eqref{def-F1}.
Then
 $$f(z)=P[\varphi](z)
+\frac{\iota}{4n}(1-|z|^2),~z\in\mathbb{B}^n,$$
 is a solution to the Poisson
equation (\ref{eq-v0-0}) with $f=\varphi$ on $\partial \mathbb{B}^n$ and $\Delta f=g_{1}=-\iota$ on $\mathbb{B}^n$, where
$$0\leq\iota<2n+\frac{4n}{1+\beta}\frac{1-\alpha}{1+\alpha}.$$
Since
\begin{align*}
\varphi(\zeta \eta)&= -F_1(-\zeta)
\\
&=\zeta\frac{\alpha+\zeta\frac{\beta+\zeta}{1+\beta\zeta}}{1+\alpha\zeta\frac{\beta+\zeta}{1+\beta\zeta}},
\end{align*}
elementary calculations yield $$D\varphi(0)\eta=\alpha~\mbox{and}~
D^2\varphi(0)(\eta^2)=2\beta(1-\alpha^{2}),$$ from which it follows that
\beqq
1+\frac{2}{1+\lambda_{\varphi,\eta}(0)}\frac{1-|D\varphi(0)\eta|}{1+|D\varphi(0)\eta|}
=1+\frac{2}{1+\beta}\frac{1-\alpha}{1+\alpha}.
\eeqq
On the other hand,
\beqq
\liminf_{r\rightarrow1^{-}}\frac{|f(\eta)-f(r\eta)|}{1-r}&=&\lim_{r\rightarrow1^{-}}
\frac{\left|P[\varphi](\eta)-P[\varphi](r\eta)-\frac{(1-r^{2})}{4n}\iota\right|}{1-r}\\
&=&\left|D\varphi(\eta)\eta-\frac{\iota}{2n}\right|\\
&=&\left|1+\frac{2}{1+\beta}\frac{1-\alpha}{1+\alpha}-\frac{\iota}{2n}\right|\\
&=&
1+\frac{2}{1+\beta}\frac{1-\alpha}{1+\alpha}-\frac{\iota}{2n}.
\eeqq
The proof of this theorem is complete.}
\qed

\subsection*{Proof of Theorem \ref{thm-chmv-1}}
It follows from \cite[Theorem 2]{Be} that all solutions to the inhomogeneous biharmonic
equation (\ref{eq-vl-1}) are given by
\beqq
f(z)&=&P[\varphi_{2}](z)+\frac{1}{2\pi}\int_{0}^{2\pi}\overline{z}e^{it}\varphi_{2}(e^{it})\frac{(1-|z|^{2})}{(1-\overline{z}e^{it})^{2}}dt\\
&&
-(1-|z|^{2})P[\varphi_{3}](z)-\frac{1}{16\pi}\int_{\mathbb{D}}g(w)G_b(z,w)dA(w),
\eeqq
where $\varphi_{3}(e^{it})=\varphi_{1}(e^{it})e^{-it}$, $dA$ denotes the Lebesgue area measure in $\mathbb{D}$
and
$$G_b(z,w)=|z-w|^{2}\log\left|\frac{1-z\overline{w}}{z-w}\right|-(1-|z|^{2})(1-|w|^{2}),
\quad z,~w\in\mathbb{D},$$
denotes the biharmonic Green function.

 Since $\varphi_{2}$ is holomorphic in $\mathbb{D}$ and continuous on $\overline{\mathbb{D}}$, we see that
$$0=\lim_{r\rightarrow1^{-}}\frac{1}{2\pi}\int_{0}^{2\pi}\overline{z}re^{it}P[\varphi_{2}](re^{it})\frac{(1-|z|^{2})}{(1-\overline{z}re^{it})^{2}}dt=
\frac{\overline{z}(1-|z|^{2})}{2\pi}\int_{0}^{2\pi}\frac{e^{it}\varphi_{2}(e^{it})}{(1-\overline{z}e^{it})^{2}}dt.$$
Consequently,
\beqq
f(z)=P[\varphi_{2}](z)
-(1-|z|^{2})P[\varphi_{3}](z)-\frac{1}{16\pi}\int_{\mathbb{D}}g(w)G_b(z,w)dA(w),
\eeqq
which, together with Theorem \ref{thm-0} (or \eqref{jm-2}), implies that

\beq\label{ju-1}
\left|f(z)\right|&\leq&\frac{|P'[\varphi_{2}](0)|+\Psi(z)}{1+\Psi(z)|P'[\varphi_{2}](0)|}|z|
+(1-|z|^{2})|P[\varphi_{3}](z)|\\ \nonumber
&&+\frac{1}{16\pi}\left|\int_{\mathbb{D}}g(w)G_b(z,w)dA(w)\right|
\eeq
where $$\Psi(z)=\frac{|z|+\frac{|P''[\varphi_{2}](0)|}{2(1-|P'[\varphi_{2}](0)|^2)}}
{1+\frac{|P''[\varphi_{2}](0)|}{2(1-|P'[\varphi_{2}](0)|^2)}|z|}|z|.$$

By \cite[Exercise 15 in Charpter 7]{G},  we see that the biharmonic Green function $G_b$ satisfies $G_b(z,w)\leq0$ for $z,w\in\mathbb{D}$.
Elementary computations lead to
\beq\label{ju-2}
\frac{1}{16\pi}\left|\int_{\mathbb{D}}g(w)G_b(z,w)dA(w)\right|&\leq&\frac{\|g\|_{\infty}}{16\pi}\int_{\mathbb{D}}\left|G_b(z,w)\right|dA(w)\\ \nonumber
&=&\frac{\|g\|_{\infty}}{16\pi}\bigg[\int_{\mathbb{D}}(1-|z|^{2})(1-|w|^{2})dA(w)\\ \nonumber
&&-
\int_{\mathbb{D}}|z-w|^{2}\log\left|\frac{1-z\overline{w}}{z-w}\right|^{2}dA(w)\bigg]\\ \nonumber
&=&\frac{\|g\|_{\infty}}{64}(1-|z|^{2})^{2},
\eeq where $\|g\|_{\infty}=\sup_{z\in\mathbb{D}}|g(z)|$.
Combining (\ref{ju-1}) and (\ref{ju-2}) gives
\beqq
|f(\eta)-f(r\eta)|
&\geq&1-\frac{|P'[\varphi_{2}](0)|+\Psi(r\eta)}{1+\Psi(r\eta)|P'[\varphi_{2}](0)|}r-(1-r^{2})|P[\varphi_{3}](r\eta)|\\
&&-\frac{\|g\|_{\infty}}{64}(1-r^{2})^{2},
\eeqq
which 
yields that
\beqq
\liminf_{r\rightarrow1^{-}}\frac{|f(\eta)-f(r\eta)|}{1-r}&\geq&
\lim_{r\rightarrow1^{-}}\frac{1-\frac{|P'[\varphi_{2}](0)|+\Psi(r\eta)}{1+\Psi(r\eta)|P'[\varphi_{2}](0)|}r}{1-r}-
\lim_{r\rightarrow1^{-}}\frac{(1-r^{2})|P[\varphi_{3}](r\eta)|}{1-r}\\
&=&1+\frac{2}{1+\frac{|P''[\varphi_{2}](0)|}{2(1-|P'[\varphi_{2}](0)|^{2})}}\frac{1-|P'[\varphi_{2}](0)|}{1+|P'[\varphi_{2}](0)|}-2|\varphi_{1}(\eta)|.
\eeqq
Next, we prove the sharpness part. For $z\in\mathbb{D}$, let $$f(z)=F_2(z)
-(1-|z|^{2})P[\varphi_{3}](z)-\frac{1}{16\pi}\int_{\mathbb{D}}g(w)G_b(z,w)dA(w),$$
which is a solution to the inhomogeneous biharmonic
equation (\ref{eq-vl-1}) with $\varphi_2=F_2$ on $\partial \mathbb{D}$.
 Here
 \begin{align*}
 F_2(z)&=e^{i\theta}\overline{\eta}z\phi_{\alpha}(-\overline{\eta}z\phi_{\beta}(-\overline{\eta}z))
 \\
 &=e^{i\theta}\overline{\eta}z
 \frac{\alpha+\overline{\eta}z\frac{\beta+\overline{\eta}z}{1+\beta\overline{\eta}z}}
 {1+\alpha\overline{\eta}z\frac{\beta+\overline{\eta}z}{1+\beta\overline{\eta}z}}
 \end{align*}
 for $z\in \mathbb{D}$,
 where $\alpha,\beta\in[0,1)$,
 \begin{align*}
 \theta&=\left\{\begin{array}{ll}
 \arg(\varphi_1(\eta)\overline{\eta})+\pi, & \mbox{if } \varphi_1(\eta)\neq 0,
 \\
 0, & \mbox{if }\varphi_1(\eta)= 0,
 \end{array}
 \right.
 \end{align*}
 \[
 2|\varphi_{1}(\eta)|<1+\frac{2}{1+\beta}\frac{1-\alpha}{1+\alpha},
 \]
 $g(z)\in C(\overline{\mathbb{D}})$.
Elementary calculations yield $$F_2'(0)=P'[\varphi_{2}](0)=e^{i\theta}\overline{\eta}\alpha~\mbox{and}~
F_2''(0)=P''[\varphi_{2}](0)=2e^{i\theta}\overline{\eta}^2\beta(1-\alpha^{2}),$$ from which it follows that
\beqq
1+\frac{2}{1+\frac{|P''[\varphi_{2}](0)|}{2(1-|P'[\varphi_{2}](0)|^{2})}}\frac{1-|P'[\varphi_{2}](0)|}{1+|P'[\varphi_{2}](0)|}
=1+\frac{2}{1+\beta}\frac{1-\alpha}{1+\alpha}.
\eeqq
On the other hand,
\beqq
\liminf_{r\rightarrow1^{-}}\frac{|f(\eta)-f(r\eta)|}{1-r}&=&\lim_{r\rightarrow1^{-}}
\frac{\left|P[\varphi_{2}](\eta)-P[\varphi_{2}](r\eta)+(1-r^{2})P[\varphi_{3}](r\eta)\right|}{1-r}\\
&=&\left|P'[\varphi_{2}](\eta)\eta+2\varphi_{1}(\eta)\overline{\eta}\right|\\
&=&\left| e^{i\theta}\left(1+\frac{2}{1+\beta}\frac{1-\alpha}{1+\alpha}\right)+2\varphi_{1}(\eta)\overline{\eta}\right|
\\
&=&
1+\frac{2}{1+\beta}\frac{1-\alpha}{1+\alpha}-2|\varphi_{1}(\eta)|.
\eeqq
The proof of this theorem is complete.
\qed

\subsection*{Proof of Theorem \ref{thm-Kalaj}}
For $\zeta\in\mathbb{D}$, let
\begin{align}
\label{def-F}
F(\zeta)=\langle f(\zeta a),b \rangle,
\end{align}
where $f:\,\Omega\to \mathbb{B}_H$ is a holomorphic mapping with $f(a)=b.$
By the assumption, we have $F(1)=1$.
Then $F$ is a holomorphic mapping of $\mathbb{D}$ into itself
such that $F$ is differentiable at $\zeta=1$ and $F(1)=1$. Elementary computations lead to
$$F'(1)=\langle Df(a)a, b\rangle\ \mbox{and}~F'(0)=\langle Df(0)a,b\rangle.$$
Let $\chi(t)=(1-t^2)+ ti$ for $t\in [-1,1]$.
Then $\chi(t)\in \mathbb{D}$ for $t\in [-1,1]\setminus \{ 0\}$
and $\chi(0)=1\in \partial \mathbb{D}$.
Therefore, ${\rm Re} \big(F(\chi(t))\big)$ attains its local maximum at $t=0$,
which implies that ${\rm Re} \big((F'(1)i)\big)=0$, where ${\rm Re}$ denotes the real part of a complex number.
Thus, $F'(1)\in \mathbb{R}$.
Moreover, we have
\begin{align*}
F'(1)={\rm Re} (F'(1))=\lim_{t\to 1^{-}}\frac{{\rm Re} (F(1))-{\rm Re}(F(t))}{1-t}\geq 0,
\end{align*}
which, together with \cite[Theorem 1.1]{Z18}, yields that
\beqq
\langle Df(a)a, b\rangle\geq
\frac{2\left| 1-\langle f(0),b \rangle\right|^2}
{1-\left|\langle f(0),b \rangle\right|^2+|\langle Df(0)a, b\rangle|}.
\eeqq

Next, assume that the equality in \eqref{eq-Schwarz-Hilbert2-new} holds.
Let $F$ be as in \eqref{def-F}.
Then $F\in \mathscr{H}(\mathbb{D}, \mathbb{D})$ with $F(1)=1$
which satisfies the assumptions of \cite[Theorem 1.1]{Z18} and the equality in \eqref{eq-Schwarz-disk}.
By \cite[Theorem 1.1]{Z18},
we deduce that
$F(\zeta)\equiv \zeta$
or
\begin{align*}
F(\zeta)\equiv \zeta\frac{\zeta+c}{1+{c}\zeta}
\end{align*}
for some $c\in [0,1)$.
Let $g(\zeta)=f(\zeta a)-F(\zeta)b$.
Then we have
$\langle b, g(\zeta)\rangle=0$ for $\zeta \in \mathbb{D}$.
Let
\begin{align*}
M(r)=\min_{|\zeta|=r}\left|  \zeta\frac{\zeta+c}{1+{c}\zeta}  \right|.
\end{align*}
Then for $\zeta$ with $|\zeta|=r\in (c,1)$,
we have
\begin{align*}
|g(\zeta)|^2=|f(\zeta a)|^2-|F(\zeta)|^2
\leq 1-M(r)^2.
\end{align*}
By the maximum principle for holomorphic functions, we have
\begin{align}
\label{g-estimate}
|g(\zeta)|^2
\leq 1-M(r)^2,
\quad
|\zeta|\leq r.
\end{align}
Let $\zeta\in \mathbb{D}$ be fixed.
Letting $r\to 1^{-}$ in \eqref{g-estimate},
we have
$g(\zeta)=0$.
Since $\zeta\in \mathbb{D}$ is arbitrary,
we obtain that
$f(\zeta)\equiv F(\zeta)b$, which implies that
$f(\zeta a)\equiv \zeta b$ or
\begin{align*}
f(\zeta a)\equiv \zeta\frac{\zeta+c}{1+{c}\zeta}b
\end{align*}
for some $c\in [0,1)$.

Finally, we prove the sharpness part.
Since $\Omega=\mathbb{B}_{X}$ is the unit ball of $X$,
for any holomorphic function $f$ of $\mathbb{D}$ into itself and for any
$l_{a}\in T(a)$ and $b\in \partial \mathbb{B}_H$,
the mapping $F(z)=f(l_{a}(z))b$ is a holomorphic mapping of $\mathbb{B}_{X}$ into $\mathbb{B}_H$.
Let $f$ be an extreme function as in \cite[Theorem 1.1]{Z18}.
Then $F(\zeta a)=f(\zeta)b$, $F(a)=b$, $\langle DF(a)a, b\rangle=f'(1)$, $\langle F(0), b\rangle=f(0)$ and $\langle DF(0)a, b\rangle=f'(0)$,
which combined with \cite[Theorem 1.1]{Z18} imply that equality in $(\ref{eq-Schwarz-Hilbert-new})$ is possible for each values of $\langle f(0),b\rangle$ and $|\langle Df(0)a, b\rangle|$
with $| \langle Df(0)a, b\rangle| \leq 1-|\langle f(0),b\rangle|^2$.
The proof of this theorem is finished.
\qed


\subsection*{Proof of Theorem \ref{thm-09-CHMV}}
It follows from Theorem {\ref{thm-h1}} that
\beqq
\|Df(z_{0})z_{0}\|_{Y}&=&\lim_{r\rightarrow 1^{-}}\left\|\frac{f(rz_{0})-f(z_{0})}{r-1}\right\|_{Y}
\geq\lim_{r\rightarrow 1^{-}}\frac{1-\|f(rz_{0})\|_{Y}}{1-r}\\
&\geq&\lim_{r\rightarrow 1^{-}}\frac{1-\frac{r+\|f(0)\|_{Y}}{1+\|f(0)\|_{Y}r}}{1-r}\\
&=&\lim_{r\rightarrow 1^{-}}\frac{1-\|f(0)\|_{Y}}{1+\|f(0)\|_{Y}r}\\
&=&\frac{1-\|f(0)\|_{Y}}{1+\|f(0)\|_{Y}}.
\eeqq

Next, we prove the sharpness part.  For  $z_{0}\in \partial\mathbb{B}_{X}$, let
 $l_{z_{0}}\in T({z_{0}})$ be a fixed element.
 For each fixed  $b\in\partial \mathbb{B}_{Y}$ and any fixed $a\in[0,1)$, define
 $$f(z)=\frac{a+l_{{z_{0}}}(z)}{1+al_{{z_{0}}}(z)}b,~z\in \mathbb{B}_{X}.$$
Then $f$ is a holomorphic mapping from $\mathbb{B}_X$ into $\mathbb{B}_Y$. 
Since $f(z_{0})=b$ and
$$f(rz_{0})=\frac{a+l_{{z_{0}}}(rz_{0})}{1+al_{{z_{0}}}(rz_{0})}b
=\frac{a+rl_{{z_{0}}}(z_{0})}{1+ral_{{z_{0}}}(z_{0})}b=\frac{a+r}{1+ra}b,$$ where $r\in[0,1)$,
we see that
\beqq
\|Df(z_{0})z_{0}\|_{Y}=\lim_{r\rightarrow 1^{-}}\left\|\frac{f(rz_{0})-f(z_{0})}{r-1}\right\|_{Y}=
\lim_{r\rightarrow 1^{-}}\frac{1-a}{1+ar}=\frac{1-\|f(0)\|_{Y}}{1+\|f(0)\|_{Y}}.
\eeqq
 The proof of this theorem is complete.
\qed

\subsection*{Proof of Theorem \ref{thm-1}}

Let $F$ be a holomorphic function from $\mathbb{D}$ into itself.
By the Schwarz-Pick lemma, we have
\be\label{eq-o-1}\left|\phi_{F(z_{0})}(F(\zeta))\right|\leq|\phi_{z_{0}}(\zeta)|,~\zeta\in\mathbb{D}.\ee
Since
\begin{align}
\label{eq-ab-estimate}
\frac{|a|-|b|}{1-|a| |b|}&\leq \left|\frac{a-b}{1-\overline{a}b}\right|,
\quad
a, b\in \mathbb{D}
\end{align}
holds (see the proof of
\cite[Proposition 2.2.2]{Kran}),
for $\zeta\in \mathbb{D}$, we have
\beqq
\frac{|F(\zeta)|-|F(z_{0})|}{1-|F(\zeta)||F(z_{0})|}\leq\left|\phi_{F(z_{0})}(F(\zeta))\right|,
\eeqq
which, together with (\ref{eq-o-1}), implies that
\be\label{eq-o-2}
|F(\zeta)|\leq\frac{|F(z_{0})|+|\phi_{z_{0}}(\zeta)|}{1+|F(z_{0})||\phi_{z_{0}}(\zeta)|}.
\ee

For any fixed $b\in\partial\mathbb{B}_Y$, let
$$F_b(\zeta)=l_{b}(f(\zeta)),~\zeta\in\mathbb{D},$$
where $l_{b}\in T(b)$. Then $F_b$ is a holomorphic function from $\mathbb{D}$ into itself with $F_b(0)=0$.
Let $$g(z)=\begin{cases}
\displaystyle \frac{F_b(z)}{z}
& \mbox{if } z\in\mathbb{D}\setminus\{0\},\\
\displaystyle F_b'(0) &\mbox{if } z=0.
\end{cases}$$
Then $g$ is also a holomorphic function on $\mathbb{D}$ into itself
or $g(\zeta)\equiv e^{i\theta}$ for some $\theta \in \mathbb{R}$.
If $g \in \mathscr{H}(\mathbb{D},\mathbb{D})$,
applying the inequality (\ref{eq-o-2}) to $g$ gives
\be\label{eq-o-3}
|F_b(\zeta)|\leq\frac{\frac{|F_b(z_{0})|}{|z_{0}|}+|\phi_{z_{0}}(\zeta)|}{1+\frac{|F_b(z_{0})|}{|z_{0}|}|\phi_{z_{0}}(\zeta)|}|\zeta|,~\zeta\in\mathbb{D}.
\ee
This inequality also holds in the case $g(\zeta)\equiv e^{i\theta}$.

Let $z$ be any fixed point in $\mathbb{D}$. Without loss of generality,
we assume that $f(z)\neq0$.
Since $$|F_b(z_{0})|=|l_{b}(f(z_{0}))|\leq\|l_{b}\|_{Y^{\ast}}\|f(z_{0})\|_{Y}=\|w_{0}\|_{Y}$$
and from the monotonic increasing property of the function $Q$  in (\ref{hyu-1}) and (\ref{eq-o-3}),
we see
that
\be\label{eq-o-4}
|l_{b}(f(z))|=|F_b(z)|\leq\frac{\frac{\|w_{0}\|_{Y}}{|z_{0}|}+|\phi_{z_{0}}(z)|}{1+\frac{\|w_{0}\|_{Y}}{|z_{0}|}|\phi_{z_{0}}(z)|}|z|.
\ee
By taking $b=f(z)/\|f(z)\|_{Y}$ in (\ref{eq-o-4}), we have
\beqq
\|f(z)\|_{Y}\leq\frac{\frac{\|w_{0}\|_{Y}}{|z_{0}|}+|\phi_{z_{0}}(z)|}{1+\frac{\|w_{0}\|_{Y}}{|z_{0}|}|\phi_{z_{0}}(z)|}|z|.
\eeqq
Next, we prove the sharpness.
If $w_0=0$,
then the mapping $f(z)=z\phi_{z_0}(z)b$,
where $b\in \partial \mathbb{B}_Y$,
attains the equality in \eqref{VM2_eq1}.
So, we consider the case $w_0\neq 0$.
Also, if $\| w_0\|_Y=|z_0|$, then
the mapping $f(z)=z w_0/z_0$
attains the equality in \eqref{VM2_eq1}.
So, we may assume that $\| w_0\|_Y<|z_0|$.
For any fixed $\alpha\in\partial\mathbb{D}$, let
	$$  f(z)=\left(z\phi_{\frac{\|w_0\|_{Y}}{z_0}}(\alpha\phi_{z_0}(z))\right)\frac{w_0}{\|w_0\|_{Y}},~z\in\mathbb{D}.
	$$
	Then $f \in \mathscr{H}(\mathbb{D},\mathbb{B}_Y)$ with $f(0)=0$, $f(z_0)=w_0$ and
	\beqq
		\|f(z)\|_{Y}&=&\big|z\phi_{\frac{\|w_0\|_{Y}}{z_0}}(\alpha\phi_{z_0}(z))\big|\\
		&=&|z| \left|\frac{\frac{\|w_0\|_{Y}}{z_{0}}-\alpha\phi_{z_0}(z)}{1-\frac{\|w_0\|_{Y}}{\overline{z}_{0}}\alpha\phi_{z_0}(z)}\right|.
	\eeqq
If we take $z=z_{1}$with $\phi_{z_0}(z_1)\neq 0$ and $\alpha=-\left(e^{i\left(\arg\frac{\|w_0\|_{Y}}{{z}_{0}}-\arg\phi_{z_0}(z_{1})\right)}\right)$, then
\beqq
		\|f(z_{1})\|_{Y}
		= \frac{\frac{\|w_0\|_{Y}}{|z_{0}|}+|\phi_{z_0}(z_{1})|}{1+\frac{\|w_0\|_{Y}}{|z_{0}|}|\phi_{z_0}(z_{1})|}|z_{1}|.
	\eeqq
The proof of this theorem is complete.
\qed

\subsection*{Proof of Theorem \ref{VM2_tthm1}}
$\phi_{f(0)} \circ f$ is a holomorphic mapping from $\mathbb{D}$ into itself which maps 0 to 0. It follows from Theorem \ref{thm-1} that
	\begin{equation}\label{VM2_eqnq}
\big|(\phi_{f(0)} \circ f)(z)\big|\leq K(z)|z|.
    \end{equation}
    It follows from \eqref{eq-ab-estimate} that, for $z\in\mathbb{D}$,
\beqq
\frac{|f(z)|-|f(0)|}{1-|f(z)||f(0)|}\leq \big|(\phi_{f(0)} \circ f)(z)\big|,
\eeqq
    which, together with (\ref{VM2_eqnq}), implies (\ref{VM2_eqnq1}).

Next, we demonstrate the sharpness of the bound.
Take $\varepsilon>0$ so that for $\alpha=(|w_0|+\varepsilon)e^{i\theta}$,
\begin{align*}
\theta&=\left\{
\begin{array}{ll}
\arg w_0, & \mbox{if } w_0\neq 0, \\
0, & \mbox{if } w_0=0,
\end{array}
\right.
\end{align*}
and $\tilde{w_0}=\phi_{\alpha}(w_0)$,
$|\tilde{w_0}|<|z_0|$ holds.
Let
\begin{align*}
f(z)&=\phi_{\alpha}\left(z\phi_{\frac{|\tilde{w_0}|}{z_0}}\left(\beta\phi_{z_0}(z)\right)\frac{\tilde{w_0}}{|\tilde{w_0}|}
\right),
\end{align*}
where $\beta \in \partial \mathbb{D}$.
Then
$f \in \mathscr{H}(\mathbb{D},\mathbb{D})$,
$f(z_0)=w_0$ and $f(0)=\alpha$.
Let $z_1\in \mathbb{D}\setminus \{ 0, z_0\}$ be such that
$\arg z_{1}=\arg z_{0}+\pi$
and let $\beta=e^{i(\pi-\arg(z_0)-\arg(\phi_{z_0}(z_1))}$.
Then
\begin{align*}
|f(z_1)|&=\left|\phi_{\alpha}\left(z_1\phi_{\frac{|\tilde{w_0}|}{z_0}}\left(\beta\phi_{z_0}(z_1)\right)\frac{\tilde{w_0}}{|\tilde{w_0}|}
\right)\right|
\\
&=\frac{K(z_1)|z_1|+|f(0)|}{1+K(z_1)|z_1||f(0)|}.
\end{align*}
 The proof of this theorem is complete.
\qed

\subsection*{Proof of Theorem \ref{VM2_thm6}}
Let $\eta=a/\|a\|_X\in\partial\mathbb{B}_X$ and let
$$F(\xi)=f(\xi\eta),~\xi\in\mathbb{D}.$$
Then $F$ is a holomorphic mapping of $\mathbb{D}$ into $\mathbb{B}_Y$
with $F(0)=0$.
Since
$
\|F(\|a\|_X)\|_Y=\|f(a)\|_{Y}\leq\delta\|a\|_X,
$
it follows from Theorem \ref{thm-1} that
\be\label{qw-1}
\left\| \frac{f(\xi \eta)}{\xi}\right\|_Y \leq\frac{|\phi_{\|a\|_{X}}(\xi)|+\delta}{1+\delta|\phi_{\|a\|_{X}}(\xi)|}.
\ee
Finally, by letting $\xi \to 0$ in (\ref{qw-1}), we obtain
$$\|Df(0)\eta\|_{Y}\leq\frac{\|a\|_{X}+\delta}{1+\delta\|a\|_{X}}.$$


	To establish the sharpness part, we proceed as follows.
	Choosing $l_{a}\in T(a)$.
Fix $b\in\partial\mathbb{B}_Y$ and $\delta\in(0,1)$.
	Let
\begin{align}
\label{eq-ex-f}
f(z)&=l_a(z)\phi_{\delta}\left(-\phi_{\| a\|_X}(l_a(z))\right)b.
\end{align}
Then $f \in \mathscr{H}(\mathbb{B}_X,\mathbb{B}_Y)$,
$f(0)=0$
and $\| f(a)\|_Y=\delta \| a\|_X$.
Since
\begin{align*}
f(\xi \eta)&=\xi \phi_{\delta}(-\phi_{\| a\|_X}(\xi))b,
\end{align*}
where $\eta=a/\| a\|_X$,
we have
\begin{align*}
Df(0)\eta&=\phi_{\delta}(-\phi_{\| a\|_X}(0))b=\frac{\|a\|_{X}+\delta}{1+\delta\|a\|_{X}}b.
\end{align*}
The proof is now complete.
\qed

\subsection*{Proof of Theorem \ref{thm-s-1}}
Let $z$ be an arbitrary fixed point in $\mathbb{B}_X\backslash\{0\}$, and assume without
loss of generality that $f(z)\neq0$.
 Let $\eta=\frac{z}{\|z\|_{X}}\in\partial\mathbb{B}_{X}$.
Define the function
$$F(\zeta)=f(\zeta\eta),~\zeta\in\mathbb{D}.$$
Then $F$ is a holomorphic mapping from $\mathbb{D}$ into $\mathbb{B}_Y$
with $F(0)=0$.
For fixed $r\in (0,1)$,
it follows from Theorem \ref{thm-1} that
\be\label{eq-o-p1}
\|F(\zeta)\|_Y\leq\frac{\frac{\|F(r)\|_Y}{r}+|\phi_{r}(\zeta)|}{1+\frac{\|F(r)\|_Y}{r}|\phi_{r}(\zeta)|}|\zeta|,~\zeta\in\mathbb{D}.
\ee
Since $\|F(r)\|_Y=|\|f(r\eta)\|_{Y},$
by  (\ref{eq-o-p1}),
we see that
\begin{align*}
\|F(\|z\|_{X})\|_Y&=\|f(z)\|_Y\leq\frac{\frac{\|f(r\eta)\|_{Y}}{r}+|\phi_{r}(\|z\|_{X})|}{1+\frac{\|f(r\eta)\|_{Y}}{r}|\phi_{r}(\|z\|_{X})|}\|z\|_{X}.
\end{align*}

 Next, we prove the sharpness.
 For any fixed 
 $r \in (0,1)$ and $b \in \partial\mathbb{B}_Y$,
 Let $f$ be as in \eqref{eq-ex-f},
 where $\| a\|_X=r$.
 Let $\eta=a/\|a\|_X\in\partial\mathbb{B}_X$.
 Then $\| f(r\eta)\|_Y=\|f(a)\|_Y=r\delta $.
 For any $t\in (0,r)$, let $z_1=t\eta$.
 Then, we have
 \begin{align*}
 \| f(z_1)\|_Y&=t\frac{\delta+\phi_{r}(t)}{1+\delta \phi_r(t)}
 =\frac{\frac{\|f(r\eta)\|_{Y}}{r}+|\phi_{r}(\| z_1\|_X)|}{1+\frac{\|f(r\eta)\|_{Y}}{r} |\phi_{r}(\| z_1\|_X)|}\| z_1\|_X.
 \end{align*}
The proof of this theorem is complete.
\qed

	



\subsection*{Proof of Corollary \ref{Vm2_thm1}}
    Let $\eta$ be any fixed point in $\partial\mathbb{B}_X$.
Let
$$F(\zeta)=f(\zeta\eta),~\zeta\in\mathbb{D}.$$
Then $F$ is a holomorphic mapping from $\mathbb{D}$ into $\mathbb{B}_Y$
with $F(0)=0$.
Applying Theorem \ref{thm-s-1} to $F$, we get
\beq\label{ppo-1}
		\frac{\|f(\zeta\eta)\|_Y}{|\zeta|} \leq \frac{\frac{\left\|f\left(\frac{r\zeta\eta}{|\zeta|}\right)\right\|_Y}{r}+|\phi_r(|\zeta|)|}{1+
\frac{\left\|f\left(\frac{r\zeta\eta}{|\zeta|}\right)\right\|_Y}{r}|\phi_r(|\zeta|)|}.
\eeq

Since $\left\|f\left(\frac{r\zeta\eta}{|\zeta|}\right)\right\|_Y\leq\omega,$
by (\ref{ppo-1}) and the monotonic increasing property of function $Q$ in (\ref{hyu-1}),
we see that
\beqq
	\frac{\| f(\zeta\eta)\|_Y}{|\zeta|} \leq \frac{\frac{\omega}{r}+|\phi_r(|\zeta|)|}{1+\frac{\omega}{r}|\phi_r(|\zeta|)|}.
\eeqq
Now letting $\zeta$ tends to $0$ gives
     \beqq
     \|Df(0)\eta\|_{Y}\leq\frac{\frac{\omega}{r}+r}{1+\omega}.
     \eeqq
 By  the arbitrariness of   $\eta$, we obtain the desired result.
\qed

\medskip
	\noindent{\bf Acknowledgments.} We sincerely appreciate  the anonymous referee for the meticulous
review and valuable suggestions, which have greatly enhanced the quality
 and readability of this paper.
The first author   was partly supported by the National Science
Foundation of China (grant no. 12571080) and  Guangxi Natural Science Foundation $\#$2026GXNSFFA00640002.
The second author was partially supported by
JSPS KAKENHI Grant Number JP22K03363.
The third author acknowledge the Department of Science and Technology (DST), Government of India, for supporting the Department of Mathematics at NIT Calicut under the FIST program, which enabled this research.
	
	\medskip

\section*{Compliance with Ethical Standards}

\noindent{\bf Conflict of Interest.} The authors declare that there is no conflict of interest regarding the publication of this paper.

\medskip

\medskip
\noindent{\bf Ethical Approval.} The authors confirm adherence to ethical standards. The submitted work is original and is not under simultaneous consideration by any other journal.

\medskip
\noindent{\bf Informed Consent.} Not applicable.

\medskip
\noindent{\bf Author Contributions.} All authors contributed equally to the conception, formulation, and preparation of this manuscript.

\normalsize

\end{document}